\documentstyle[12pt,twoside]{article}
\newtheorem{theorem}{\bf Theorem}[section]
\newtheorem{proposition}[theorem]{\bf Proposition}

\newtheorem{corollary}[theorem]{\bf Corollary}

\input{amssym}

\date{}
\begin{document}

\title{{\Large\bf Cohomological properties of different types of weak amenability}}

\author{{\normalsize\sc M. J. Mehdipour and A. Rejali\footnote{Corresponding author}}}
\maketitle

{\footnotesize  {\bf Abstract.} In this paper, we deal with cohomological properties of weak amenability, cyclic amenability, cyclic weak amenability and point amenability of Banach algebras. We look at some hereditary properties of them and show that continuous homomorphisms with dense range preserve cyclically weak amenability, however, weak amenability and cyclically amenability are preserved under certain conditions. We also study these cohomological properties of the $\theta-$Lau product  $A\times_\theta B$ and the projective tensor product $A\hat{\otimes} B$. Finally, we investigate the cohomological properties of $A^{**}$ and establish that  cyclically weak amenability of $A^{**}$ implies cyclically weak amenability of $A$. This result is true for point amenability instead of cyclically weak amenability.}
                     %-----------------------------------------
                     %-----------------------------------------
                     %-----------------------------------------
{}
                     %-----------------------------------------
                     %-----------------------------------------
                     %-----------------------------------------
                    % \tableofcontents

\section{\normalsize\bf Introduction}

For a Banach algebra $A$, a linear operator $D: A\rightarrow A^*$ is said to be a \emph{derivation} if $D(ab)=D(a)\cdot b + a\cdot D(b)$ for all $a, b\in A$. If for every $a, b\in A$  
$$
\langle D(a), b\rangle+\langle D(b), a\rangle=0,
$$
$D$ is called \emph{cyclic}. One can prove that $D$ is a cyclic derivation if and only if $\langle D(a), a\rangle=0$ for all $a\in A$. In the case where $A$ is unital, then  $D$ is a cyclic derivation if and only if $\langle D(a), 1_A\rangle=0$ for all $a\in A$; see \cite{mr4}.
For every $F\in A^*$, the bounded linear operator $a\mapsto F\cdot a-a\cdot F$ is called an \emph{inner derivation} and is denoted by $\hbox{ad}_F$, where $\langle F\cdot a, x\rangle=\langle F, ax\rangle$ and
$\langle a\cdot F, x\rangle=\langle F, xa\rangle
$
for all $x\in A$.

A Banach algebra $A$ is called \emph{weakly amenable} if inner derivations are the only continuous derivation from $A$ into $A^*$. 
Also, $A$ is called \emph{ cyclically weakly amenable} if every continuous derivation from $A$ into $A^*$ is cyclic, and $A$ is said to be \emph{cyclically amenable} if every cyclic derivation from $A$ into $A^*$ is inner.  

Let us recall that a linear functional $d$ on $A$ is called a \emph{continuous point derivation of} $A$ \emph{at} $\varphi\in\Delta_0(A)$ if for every $a, b\in A$, we have
$$
d(ab)=d(a)\varphi(b)+\varphi(a) d(b),
$$
where $\Delta_0(A):=\Delta(A)\cup\{0\}$ and $\Delta(A)$ is the character space of $A$.
We call  $A$ \emph{point amenable} (0-\emph{point amenable}) if there is no non-zero continuous point derivation of $A$ at $\Delta(A)$ (respectively, $\Delta_0(A)$). 

 Cohomology properties of Banach algebras studied by many authors \cite{bcd, cl, f, gro, gro1, j1, m}. For example, Bade, Curtis and Dales \cite{bcd} introduced the concept of weak
amenability for commutative Banach algebras. Johnson [25] gave this concept for arbitrary
Banach algebras. Gronbaek \cite{gro} studied the hereditary properties of weak amenability. He \cite{gro1} introduced and investigated the concept of cyclically amenablity. Recently, the authors \cite{mr4} introduced the notions of cyclically weak amenability and point amenability for  Banach algebras. They proved that $A$ is weakly amenable if and only if $A$ is both cyclically weakly amenable and cyclically amenable. They showed that if $A$ is commutative, then weakly amenability and cyclically weak amenability of $A$ coincide. They also proved that every Banach algebra with non-empty character space is cyclically weakly amenable if and only if $A$ is 0-point amenable; or equivalently, it is point amenable and essential. In \cite{mr4}, it is shown that weak amenability, cyclically weak amenability and point amenability of commutative unital Banach algebras are equivalent.

In this paper we investigate the hereditary properties of cyclically weak amenability and point amenability. We also find some properties concerning weak amenability and cyclically amenability. In Section 2, we  study some hereditary properties. We verify that continuous homomorphisms with dense range are cyclically weak amenability  preserving, but they are not neither weak amenability nor cyclic amenability preserving in general. We give condition under which continuous homomorphisms with dense range are weak amenability and cyclic amenability preserving. In Section 3, we study weak amenability, cyclic amenability, cyclic weak amenability and point amenability of Banach algebras $A\times_\theta B$, $A\oplus B$,  $A^\sharp$ and $A\hat{\otimes}B$. In Section 4, we prove that point amenability of $A^{**}$ implies point amenability of $A$, however, this result is true for cyclically weakly amenable and cyclically amenable under certain conditions.

\section{\normalsize\bf Homomorphism and retraction maps on Banach algebras}

Let $A_1$ and $A_2$ be Banach algebras and $\phi: A_1\rightarrow A_2$ be a continuous homomorphism with dense range. Then $\phi$ is said to be \emph{\emph{(}respectively,  weak, cyclic, cyclically weak, point, 0-point\emph{)} amenability preserving} if (respectively,  weak, cyclic, cyclically weak, point, 0-point) amenability of $A_1$ implies  (respectively, cyclic, cyclically weak, point, weak) amenability of $A_2$. Some authors investigated these properties. For instance, Paterson \cite{pat} proved that $\phi$ is always amenability preserving. Groenbaek studied weak amenability preserving. He \cite{gro} showed that if $A_1$ and $A_2$ are commutative, then $\phi$ is weak amenability preserving. In the case where $A_1$ and $A_2$ are  noncommutative, he \cite{gro2} gave some sufficient conditions. Lau and Loy \cite{ll} showed that $\phi$ is point amenability preserving.

\begin{theorem}\label{mm} Let $A_1$ and $A_2$ be Banach algebras and $\phi: A_1\rightarrow A_2$ be a continuous homomorphism with dense range. The following statements hold.

\emph{(i)} $\phi$ is cyclically weak amenability preserving.

\emph{(ii)} If $A_2^*\circ\phi=A_1^*$, then $\phi$ is cyclic amenability preserving.

\emph{(iii)} If $A_2^*\circ\phi=A_1^*$, then $\phi$ is weak amenability preserving.

\emph{(iv)} $\phi$ is 0-point amenability preserving.
%\emph{(i)} If $D: A_2\rightarrow A_2^*$ is a continuous derivation, then $\phi^*\circ D\circ\phi$ is a continuous derivation from $A_1$ into $A_1^*$.

%\emph{(ii)}  If $A_1$ is a weakly amenable A_2anach algebra, then %there exists $\Phi\in A_1^*$ such that $\langle D(\phi(x), \phi(y)\rangle= \langle \hbox{ad}_\Phi(x), y\rangle$ for all $x, y\in A_1$.
%$\phi^*\circ D\circ\phi$ is inner.

%\emph{(iii)} If $A_1$ is a weakly amenable A_2anach algebra and every inner derivation from $A_1$ into $A_1^*$ is determined by $Q\circ\phi$, then $A_2$ is weakly amenable, where $Q\in A_2^*$.

%\emph{(ii)}
%\emph{(ii)} If $A_1$ is a weakly amenable A_2anach algebra, then $A_2$ is weakly amenable.

%\emph{(iii)} If $A_1$ is weakly amenable, then every closed subalgebra of $A_1$ is weakly amenable.
\end{theorem}
{\it Proof.} Let $D: A_2\rightarrow A_2^*$ be a continuous derivation. In view of Lemma 2.2 in \cite{ll}, $\tilde{D}:=\phi^*\circ D\circ\phi$ is a continuous derivation from $A_1$ into $A_1^*$.

(i) Assume that $A_1$ is  cyclically weakly amenable. Then $\tilde{D}$ is cyclic. So for every $a_1, x_1\in A_1$
$$
\langle \tilde{D}( a_1), x_1\rangle+\langle \tilde{D}(x_1), a_1\rangle=0.
$$
Hence
$$
\langle D(\phi(a_1)), \phi(x_1)\rangle+\langle D(\phi( x_1))\phi(a_1)\rangle=0.
$$
Since $\phi$ is continuous homomorphism with dense range, it follows that
$$
\langle D(a_2), x_2\rangle+\langle D(x_2), a_2\rangle=0
$$
for all $a_2, x_2\in A_2$. That is, $D$ is cyclic. That is, (i) holds.

(ii) Let $A_1$ be cyclic amenable and $A_2^*\circ\phi=A_1^*$. It is easy to see that if $D$ is cyclic. Then $\tilde{D}$ is cyclic. Thus $\tilde{D}$ is inner and so there exists a functional $F\in A_2^*$ such that $\tilde{D}=\hbox{ad}_{F\circ\phi}$. Hence
for every $a_1, x_1\in A_1$, we have
\begin{eqnarray*}
\langle D(\phi(a_1)), \phi(x_1)\rangle&=&F\circ\phi(a_1x_1)-F\circ\phi(x_1a_1)\\
&=&F(\phi(a_1)\phi( x_1))- F(\phi( x_1) \phi(a_1)).
\end{eqnarray*}
This implies that
$$
\langle D(a_2), x_2\rangle=F(a_2x_2)-F(x_2a_2)
$$
for all $a_2, x_2\in A_2$. Hence $A_2$ is weakly amenable. That is (ii) holds.

(iii) The statements  (i) and (ii) together with Theorem 4.1 in \cite{mr4}  prove (iii).

(iv) An argument like that used in the proof of Lemma 2.1 in \cite{ll} proves (iv).
$\hfill\square$\\

As an immediate consequence of  Lemma 2.1 in \cite{ll} and Theorem \ref{mm} we give the following result.

\begin{corollary} Let $A$ be a Banach algebra and $I$ be a closed ideal in $A$. Then the following statements hold.

\emph{(i)} If $A$ is  cyclically weakly amenable, then $A/I$ is  cyclically weakly amenable.

\emph{(ii) }If $A$ is point amenable, then $A/I$ is point amenable.
\end{corollary}

\begin{corollary}\label{havig} Let $A$ be a Banach algebra and $B$ be a dense subalgebra of $A$. If $B$ is weakly amenable, then $A$ is weakly amenable. This statement is true if weakly amenable is replaced by cyclically amenable, cyclically weakly amenable or point amenable.
\end{corollary}
{\it Proof.} Let $\ell: B\rightarrow A$ be the inclusion map. Then $\ell$ is a continuous epimorphism with dense range. If $g\in B^*$, then there exists $f\in A^*$ such that $f|_B=g$. So $f\circ\ell=g$. It follows that $A^*\circ\ell=B^*$. These facts together with Theorem \ref{mm} prove the result.$\hfill\square$\\

Let $(A, \|.\|_A)$ and $(B, \|.\|_B)$ be Banach algebras. Then $B$ is called an \emph{abstract Segal algebra} if $B$ is a dense left ideal in $A$ and there exist positive numbers $m$ and $n$ such that
$$
\|b\|_A\leq m\|b\|_B\quad\hbox{and}\quad\|ab\|_B\leq n\|a\|_A\|b\|_B.
$$
for all $a, b\in B$.

\begin{corollary} Let $B$ be an abstract Segal algebra with respect to  Banach algebra $A$. If $A$ is weakly amenable, then so does $B$. This statement is true for cyclically amenable, cyclically weakly amenable or point amenable instead of wealy amenable.
\end{corollary}
%{\it Proof.} It is proved in \cite{ann} that $\Delta(A)=\Delta(B)$. Since every continuous point derivation of $B$ can
%be extended to a continuous point derivation of $A$. Hence $B$ is point amenable.$\hfill\square$\\

For  a Banach algebra $A$, a continuous linear functional $p\in (A\hat{\otimes}A)^*$  is called a \emph{quasi-additive functional of} $A$ if  for every $a, b, x\in A$,
$$
p(ax\otimes b)=p(a\otimes xb)+p(x\otimes ba).
$$
Furthermore, $p$  is called \emph{inner} if there exists $F\in A^*$ such that $$p(a\otimes b)=F(ab-ba)$$ for all $a, b\in A$. A quasi-additive functional $p$ of $A$ is called \emph{cyclic} if $p(a\otimes a)=0$ for all $a\in A$. One can prove that $A$ is weakly amenable if and only if every quasi-additve functional of $A$ is inner. Also, $A$ is cyclically amenable if and only if every cyclic quasi-additive functional of $A$ is inner; see \cite{mr4}.

Let $A_1$ and $A_2$ be Banach algebras. A continuous homomorphism $\phi: A_1\rightarrow A_2$ is called a \emph{retraction} if  there exists a continuous homomorphism $\psi: A_2\rightarrow A_1$ such that $\phi\circ\psi=I_{A_2}$.

\begin{theorem}\label{react}  Let $\phi$ be a retraction from Banach algebra $A_1$ into $A_2$. Then the following statements hold.

\emph{(i)} If $A_1$ is cyclic amenable, then $A_2$ is cyclic amenable.

\emph{(ii)} If $A_1$ is weakly amenable, then $A_2$ is weakly amenable.
\end{theorem}
{\it Proof.} Let $p$ be a cyclic quasi-additive functional of $A_2$. Then $p\circ (\phi\otimes\phi)$ is a cyclic quasi-additive functional of $A_1$. Thus there exists $F\in A_1^*$ such that 
$$
p(\phi(a_1)\otimes\phi(x_1))=F(a_1x_1-x_1a_1)
$$
for all $a_1, x_1\in A_1$. Since $\phi$ is a retraction, $\phi\circ\psi=I_{A_2}$ for some continuous homomorphism $\psi$ from $A_2$ into $A_1$. For every $a_2, x_2\in A_2$, we have
\begin{eqnarray*}
p(a_2\otimes x_2)&=&p(\phi\psi(a_2)\otimes \phi\psi(x_2))\\
&=&p\circ(\phi\otimes\phi)(\psi(a_2)\otimes\psi(x_2))\\
&=&F(\psi(a_2)\psi(x_2)-\psi(x_2)\psi(a_2))\\
&=&F\circ\psi(a_2x_2-x_2a_2).
\end{eqnarray*}
This shows that $p$ is inner. So (i) holds. The authors \cite{mr4} proved that a Banach algebra is weakly amenable if and only if it is cyclically amenable and cyclically weakly amenable. This facts together with part (i) and Theorem \ref{mm} (i) proves (ii). $\hfill\square$

\section{\normalsize\bf Cohomological properties of Lau Product of Banach algebras}

 Let $A_1$ and $A_2$ be Banach algebras and $\theta\in\Delta_0(A_2)$. Let us recall that the $\theta-$Lau product $A_1\times_\theta A_2$ is the Cartesian product
 $A_1\times A_2$ with the product 
$$
(a_1, a_2)(x_1, x_2)=(a_1x_1+ a_1\theta(x_2)+x_1\theta(a_2), a_2x_2)
$$
and the norm $\|(a_1, x_1)\|=\|a_1\|+\|x_1\|$. Let us remark that if $\theta=0$, then $A_1\times_\theta A_2$ is the usual direct product $A_1\oplus A_2$ and if $A_2=\Bbb{C}$ and $\theta$ is the identity map on $\Bbb{C}$, then $A_1\times_\theta\Bbb{C}$ is $A_1^\sharp$, the unitization of $A_1$.

\begin{theorem}\label{sandali124} Let $A_1$ and $A_2$ be Banach algebras and $\theta\in\Delta_0(A_2)$. Then the following statements hold.

\emph{(i)} $A_1\times_\theta A_2$ is 0-point amenable if and only if  $A_1$ and $A_2$ are 0-point amenable.

\emph{(ii)} $A_1\oplus A_2$ is 0-point amenable if and only if $A_1$ and $A_2$ are 0-point amenable.

\emph{(iii)} $A_1$ is 0-point amenable if and only if $A_1^\sharp$ is 0-point amenable.
\end{theorem}
{\it Proof.}  Let $d: A_1\rightarrow{\Bbb C}$ be a continuous point derivation at $\varphi\in\Delta_0(A_1)$. Define the bounded linear functionals $\bar{d}$ and $\bar{\varphi}$ on $A_1\times_\theta A_2$ by 
$$
\bar{d}(a_1, a_2)=d(a_1)\quad\hbox{and}\quad\bar{\varphi}(a_1, a_2)=\varphi(a_1)+\theta(a_2).
$$
It is routine to check that $\bar{d}$ is a continuous point derivation at $\bar{\varphi}\in\Delta_0(A_1\times_\theta A_2)$. So if $A_1\times_\theta A_2$ is 0-point amenable, then $\bar{d}=0$ and hence $d=0$. This shows that $A_1$ is 0-point amenable. To complete the proof, we only note that the mapping $(a_1, a_2)\mapsto a_2$ is a continuous epimorphism from $A_1\times_\theta A_2$ onto $A_2$ and it is 0-point amenable preserving. Thus $A_2$ is 0-point amenable.

Conversely, let $d\in(A_1\times_\theta A_2)^*$ be a continuous point derivation at $\varphi\in\Delta_0(A_1\times_\theta A_2)$. Put $\tilde{d}_i:=d|_{A_i}$ and $\tilde{\varphi}_i=\varphi|_{A_i}$, for $i=1, 2$. Then $\tilde{d}_i$ is a continuous point derivation at $\tilde{\varphi}_i\in \Delta_0(A_i)$. Since $A_i$ is a 0-point amenable, $\tilde{d}_i=0$ for $i=1, 2$. It follows that $d=0$. Consequently, $A_1\times_\theta A_2$ is 0-point amenable. Therefore, (i) holds.  Since 
$$
A_1\oplus A_2=A_1\times_0 A_2\quad\hbox{and}\quad A_1^\sharp=A_1\times_{id}{\Bbb C},
$$
the statements (ii) and (iii) follow from (i). $\hfill\square$\\

From Theorem 3.3 in \cite{mr4} and Theorem \ref{sandali124} we infer the following result.

\begin{corollary}\label{hadi} Let $A_i$  be a Banach algebra with $\Delta(A_i)\neq\emptyset$, for $i=1, 2$. If $\theta\in\Delta(A_2)$, then the following statements hold.

\emph{(i)} $A_1\times_\theta A_2$ is cyclically weakly amenable if and only if $A_1$ and $A_2$ are cyclically weakly amenable.

\emph{(ii)} $A_1$ is cyclically weakly amenable if and only if $A_1^\sharp$ is cyclically weakly amenable.
\end{corollary}

\begin{theorem} Let $A_1$ and $A_2$ be  Banach algebras and $\theta\in\Delta_0(A_2)$. Then the following statements hold.

\emph{(i)} If $A_1\times_\theta A_2$ is point amenable, then $A_1$ and $A_2$ are point amenable.

\emph{(ii)} If $A_1\oplus A_2$ is point amenable, then $A_1$ and $A_2$ are point amenable.

\emph{(iii)} If $A_1^\sharp$ is point amenable, then $A_1$ is point amenable.\\
The converse of above statements holds when $A_1$ and $A_2$ are essential.
\end{theorem} 
{\it Proof.}  A similar argument to that in Theorem \ref{sandali124} (i) shows that point amenability $A_1\times_\theta A_2$ implies point amenability $A_1$ and $A_2$. The converse
follows from Theorem \ref{sandali124} and the fact that an essential Banach algebra is point amenable if and only if it is 0-point amenable.$\hfill\square$

\begin{theorem}\label{tabestan} Let $A_1$ and $A_2$ be Banach algebras and $\theta\in\Delta(A)$. Then the following statements hold.

\emph{(i)} If $A_1\times_\theta A_2$ is cyclic  amenable, then  $A_2$ is cyclic amenable.

\emph{(ii)} $A_1$ is cyclically amenable if and only if $A_1^\sharp$ is cyclically amenable.
\end{theorem}
{\it Proof.} For $i= 1, 2$, let $\pi_i: A_1\times A_2\rightarrow A_i$ and $\iota_i: A_i\rightarrow A_1\times A_2$ be the canonical projection and the canonical injection maps, respectively. Then $\pi_i\circ\iota_i=I_{A_i}$. Thus $\pi_i$ is reaction and by Theorem \ref{react} the the statement (i) holds.

For (ii), let $\tilde{D}: A_1^\sharp\rightarrow (A_1^\sharp)^*=A_1^*\oplus\Bbb{C}$  be a  continuous cyclic derivation. Define the continuous cyclic derivation  $U: A_1\rightarrow A_1^*$ by
$$
\langle U(a), b\rangle=\langle\tilde{D}(a, 0), (b, 0)\rangle.
$$
If $A_1$ is cyclically amenable, then  $U=\hbox{ad}_F$ for some $F\in A_1^*$. Define the bounded linear functional $\tilde{F}$ on $A_1^\sharp$ by $\tilde{F}(a, r)=F(a)$ for all $a\in A_1$ and $r\in\Bbb{C}$. Thus $\tilde{D}=\hbox{ad}_{\tilde{F}}$, because, for every $a, b\in A_1$ and $r, s\in\Bbb{C}$, we have
\begin{eqnarray*}
\langle\hbox{ad}_{\tilde{F}}(a,r), (b,s)\rangle&=&\langle \tilde{F}, (a,r)(b,s)-(b,s)(a,r)\rangle\\
&=&\langle\tilde{F}, ( ab+ rb+ sa)-
( ba+ rb+sa)\rangle\\
&=&\langle F, ab-ba\rangle\\
&=&
\langle\hbox{ad}_F(a), b\rangle\\
&=&\langle U(a), b\rangle.
\end{eqnarray*}
Since $\tilde{D}$ is cyclic, we have
\begin{eqnarray*}
\langle \tilde{D}(a,r), (b,s)\rangle&=&
\langle \tilde{D}(a,r), (b,0)\rangle+\langle \tilde{D}(a,r), (0,s)\rangle\\
&=&
\langle \tilde{D}(a,0), (b,0)\rangle+\langle \tilde{D}(0,r), (b,0)\rangle\\
&+&\langle \tilde{D}(a,0), (0,s)\rangle+\langle\tilde{D}(0,r), (0,s)\rangle\\
&=&\langle U(a), b\rangle.
\end{eqnarray*}
Therefore, $A_1^\sharp$ is cyclically amenable. The converse follows from (i).$\hfill\square$\\

Sanjani Monfared proved that if   $A_1\times_\theta A_2$ is  weakly amenable, then $A_1$ is weakly amenable and $A_2$ is cyclically amenable; see Theorem 2.11 in \cite{mon}. In the next result, we give an improvement of it.

\begin{theorem}\label{min}  Let $A_1$ and $A_2$ be Banach algebras and $\theta\in\Delta(A_2)$. Then the following statements hold.

\emph{(i)} If  $A_1\times_\theta A_2$ is  weakly amenable, then $A_1$ and $A_2$  are weakly amenable.

\emph{(ii)} If $A_i$ is weakly amenable with $\Delta(A_i)\neq\emptyset$, for $i=1, 2$, then $A_1\times_\theta A_2$ is  cyclically weakly amenable.

\emph{(iii)} Assume that $A_i$ is commutative and $\Delta(A_i)\neq\emptyset$, for $i=1, 2$. Then $A_1\times_\theta A_2$ is  weakly amenable if and only if $A_1$ and $A_2$ are weakly amenable.
\end{theorem}
{\it Proof.} (i) Assume that $A_1\times_\theta A_2$ is  weakly amenable. Then $A_1\times_\theta A_2$ is  both cyclically amenable and cyclically weakly amenable. From  Corollary \ref{hadi} and Theorem \ref{tabestan} we conclude that that $A_2$  is  both cyclically amenable and cyclically weakly amenable. Thus $A_2$ is weakly amenable. From this and Theorem 2.11 in \cite{mon} we see that (i) holds.

(ii) If  $A_1$ and $A_2$ are weakly amenable, then they are cyclically weakly amenable. By Corollary \ref{hadi}, $A_1\times_\theta A_2$ is  cyclically weakly amenable. 

(iii) This follows from (i) and (ii).$\hfill\square$\\

For a Banach algebra $A$, it is proved that if $A^\sharp$ is weakly amenable and $H^2(A, \Bbb{C}_0)=\{0\}$, then $A$ is weakly amenable; see Proposition 2.8.67 in \cite{d}. In the case where $A$ is commutative, it is shown that $A^\sharp$ is weakly amenable if and only if $A$ is weakly amenable; see Corollary 2.8.70 in \cite{d}. Gronback \cite{gro2} established that if $A$ has a bounded approximate identity, then weak amenability of $A^\sharp$ implies that $A$ is weakly amenable.

\begin{corollary}  Let $A$ be a Banach algebra.  If $A$ is weakly amenable, then $A^\sharp$ is weakly amenable. The converse holds if $\Delta(A)\neq\emptyset$.
\end{corollary}
{\it Proof.} In view of  Theorem \ref{min},  weak amenability of $A$ implies weak amenability of $A^\sharp$. To complete the proof, let  $A^\sharp$ be weakly amenable. Then $A^\sharp$ is both cyclically amenable and cyclically weakly amenable.  So if  $\Delta(A)\neq\emptyset$, then Corollary \ref{hadi}  shows that $A$ is weakly amenable.$\hfill\square$

\begin{corollary}  Let $A$ be a semisimle Banach algebra.  Then $A$ is weakly amenable if and only if $A^\sharp$ is weakly amenable. 
\end{corollary}

Lau and Loy \cite{ll} proved that if Banach algebra $A$ be the direct summand of a closed
subalgebra $B$ and closed two-sided ideal $I$, then the weak amenability of $A$ implies the weak amenability of $B$; see also Proposition 2.4 in \cite{gro2}. We now give an improvment of this result.

\begin{theorem}\label{su123}  Let $A_1$ and $A_2$ be Banach algebras. Then the following statements hold.

 \emph{(i)} $A_1\oplus A_2$ is  cyclically weakly amenable if and only if $A_1$ and $A_2$ are  cyclically weakly amenable.

\emph{(ii)} $A_1\oplus A_2$ is  cyclically amenable if and only if $A_1$ and $A_2$ are cyclically amenable.

 \emph{(iii)} $A_1\oplus A_2$ is weakly amenable if and only if $A_1$ and $A_2$ are weakly amenable.
%%\emph{(v)} $A_1\oplus A_2$ is 0-point amenable if and only if $A_1$ and $A_2$ are 0-point amenable.
\end{theorem}
{\it Proof.} (i) Let $D: A_1\oplus A_2\rightarrow (A_1\oplus A_2)^*$ be a continuous derivation. Define the continuous derivations $D_1: A_1\rightarrow A_1^*$ and $D_2: A_2\rightarrow A_2^*$ by
$$
\langle D_1(a_1), x_1\rangle=\langle D(a_1, 0), (x_1, 0)\rangle\quad\hbox{and}\quad\langle D_2(a_2), x_2\rangle=\langle D(0, a_2), (0, x_2)\rangle.
$$
Then $D_1$ and $D_2$ are cyclic. Note the mapping $F\mapsto(F_1, F_2) $ is an isometric  isomorphism from  $ (A_1\oplus A_2)^*$ onto $A_1^*\oplus A_2$, where 
$$
F_1(a_1)=F(a_1, 0)\quad\hbox{and}\quad F_2(a_2)=F(0, a_2)
$$
for all $a_1\in A_1$ and $a_2\in A_2$. This implies that 
$$
\langle D(a_1, 0), ( 0, x_2)\rangle=\langle D(0, a_2), (x_1, 0)\rangle=0
$$
for all $a_1\in A_1$ and $a_2\in A_2$.
 Hence
\begin{eqnarray*}
\langle D( a_1,a_2),( x_1 ,x_2 )\rangle&+&\langle D(x_1 ,x_2 ), (a_1,a_2)\rangle\\&=&
\langle D(a_1,a_2), (x_1 ,0)\rangle+\langle D(a_1,a_2), (0,x_2 )\rangle\\
&+&
\langle D(x_1 ,x_2 ), (a_1,0)\rangle+\langle D(x_1 ,x_2 ), (0,a_2)\rangle\\
&=&
\langle D(a_1,0), (x_1 ,0)\rangle+\langle D(0,a_2), ( x_1  ,0)\rangle\\
&+&
\langle D(a_1,0), (0,x_2 )\rangle+\langle D(0,a_2), (0,x_2 )\rangle\\
&+&
\langle D(x_1 ,0), (a_1,0)\rangle+\langle D(0,x_2 ), (a_1,0)\rangle\\
&+&
\langle D(x_1 ,0), (0,a_2)\rangle+\langle D(0,x_2 ), (0,a_2)\rangle\\
&=&
\langle D_1(a_1), x_1 \rangle+\langle D_1(x_1 ), a_1\rangle\\
&+&
\langle D_2( a_2), x_2 \rangle+\langle D_2(x_2 ), a_2\rangle\\
&=&0
\end{eqnarray*}
for all $a_1, x_1\in A_1$ and $a_2, x_2\in A_2$.
This shows that $D$ is cyclic. Therefore, $A_1\oplus A_2$ is cyclically weakly amenable.

 Let $\tilde {D}_1: A_1\rightarrow A_1^*$ be a continuous derivation.  We define the linear operator $\tilde{D}: A_1\oplus A_2\rightarrow(A_1\oplus A_2)^*$ by
$$
\langle \tilde{D}(a_1, a_2), (x_1, x_2)\rangle=\langle \tilde{D_1}(a_1), x_1\rangle
$$
for all $a_1, x_1\in A_1$ and $a_2, x_2\in A_2$. Then $\tilde{D}$ is a continuous derivation.
If $A_1\oplus A_2$ be  cyclically weakly amenable, then $\tilde{D}$ is cyclic. So for every $a_1, x_1\in A_1$ we have
\begin{eqnarray*}
\langle \tilde{D_1}(a_1), x_1\rangle&+&\langle \tilde{D_1}(x_1), a_1\rangle\\
&=&\langle \tilde{D}(a_1, 0), (x_1, 0)\rangle\\
&+&\langle \tilde{D}(x_1, 0), (a_1, 0)\rangle=0.
\end{eqnarray*}
This shows that $\tilde{D_1}$ is cyclic. Hence $A_1$ is  cyclically weakly amenable. Similarly, one can prove that $A_2$ is  cyclically weakly amenable.

(ii) Let $A_1$ and $A_2$ are cyclically amenable. Assume that
$D: A_1\oplus A_2\rightarrow A_1^*\oplus A_2^*$ is a cyclic derivation. Define $D_1: A_1\rightarrow A_1^*$ and $D_2: A_2\rightarrow A_2^*$ by
$$
\langle D_1(x), a\rangle=\langle D(x, 0), (a, 0)\rangle\quad\hbox{and}\quad\langle D_2(y), b\rangle=\langle D(0, y), (0, b)\rangle.
$$
Then
$$
\langle D_1(x), a\rangle+\langle D_1(a), x\rangle=\langle D(x, 0), (a, 0)\rangle+\langle D(a, 0), (x, 0)\rangle=0
$$
for all $a, x\in A_1$. So $D_1$ is a cyclic derivation. Similarly, $D_2$ is a cyclic derivation. Hence $D_1$ and $D_2$ are inner. Thus $D_1=\hbox{ad}_{F_1}$ and $D_2=\hbox{ad}_{F_2}$ for some $F_1\in A_1^*$ and $F_2\in A_2^*$. If $F(x, y)=F_1(x)+F_2(y)$, then $F\in (A_1\oplus A_2)^*$ and $D=\hbox{ad}_F$. That is, $A_1\oplus A_2$ is cyclically amenable. Theorem \ref{tabestan} (i) proves the converse.
%%(v) Let $A_1\oplus A_2$ be  point amenable and let $d_1$ be a continuous point derivation of $A_1$ at %%$\varphi_1\in\Delta(A_1)\cup\{0\}$. Define the linear functionals $d,\varphi: A_1\oplus A_2\rightarrow {\Bbb C}$ by $d(a_1, %%a_2)=d_1(a_1)$ and $\varphi(a_1, a_2)=\varphi_1(a_1)$. It is obvious that $d$ is a continuous point derivation of $A_1\oplus %%A_2$ at $\varphi$. Hence $d=0$ and so $d_1=0$. This shows that $A_1$ is point amenable. Similarly, $A_2$ is point amenable.
%%Conversely, let $A_1$ and $A_2$ be point amenable. Assume that $d$ be a continuous point amenable at %%$\varphi\in\Delta(A_1\oplus A_2)\cup\{0\}$. Define the linear functionals $d_1$ and $\phi_1$ on $A_1$ by
%%$$
%%d_1(a_1)=d(a_1, 0)\quad\hbox{and}\quad\varphi_1(a_1)=\varphi(a_1, 0)
%%$$
%%for all $a_1\in A_1$. Then $d_1$ is a continuous point derivation of $A_1$ at $\varphi_1$. Since $A_1$ is point amenable, %%$d_1=0$. A similar argument applies to the linear functionals $d_2, \varphi_2: A_2\rightarrow{\Bbb C}$ defined by
%%$$
%d_2(a_2)=d(0, a_2)\quad\hbox{and}\quad\varphi_2(a_2)=\varphi(0, a_2)
%$$
%to conclude that $d_2=0$. But,
%$$
%d(a_1, a_2)=d_1(a_1)+d_2(a_2)=0
%$$
%for all $a_1\in A_1$ and $a_2\in A_2$. That is, $d=0$.

(iii) This follows from Theorem 4.1 in \cite{mr4} and the statements (ii) and (iii).
$\hfill\square$\\

Let $A_i$ be a Banach algebra, for $i=1, 2, 3$.  Let us recall that a short exact sequence 
\begin{eqnarray}\label{soheilag}
0\rightarrow A_1\rightarrow A_2\rightarrow A_3\rightarrow 0
\end{eqnarray}
is called \emph{split} if there exists a homomorphism $h: A_3\rightarrow A_2$ such that $h\circ g=I_{A_2}$. It is easy to see that the sequence (\ref{soheilag}) is split if and only if there exists a homomorphism $k:A_2\rightarrow A_1$ such that $f\circ k=I_{A_2}$; or equivalently, $A_2=A_1\oplus A_3$. One can prove that if $A_1$ is injective or $A_3$ is projective, then the sequence (\ref{soheilag}) is split. These facts together with Theorem \ref{su123} prove the following result.

\begin{corollary} \label{sail1} Let $0\rightarrow A_1\rightarrow A_2\rightarrow A_3\rightarrow 0$ be a short exact sequence of Banach algebras, and let $A_2$ be weakly amenable. Then the following statements hold.

\emph{(i)} If the given sequence is split, then $A_1$ and $A_3$ are weakly amenable.

\emph{(ii)} If $A_1$ is injective or $A_3$ is projective, then $A_1$ and $A_3$ are weakly amenable.
\end{corollary}

\begin{corollary}\label{sail2} Let $A$ be a weakly amenable Banach algebra, and let $I$ be an ideal of $A$.  If $I$ is injective or $A/I$ is projective, then $I$ and $A/I$ are weakly amenable.
\end{corollary}
{\it Proof.} It suffices to note that the sequence
$$
0\rightarrow I\stackrel{\iota}\rightarrow A\stackrel{\pi}\rightarrow A/I\rightarrow 0
$$
is a short exact sequence, where $\iota$ is the inclusion map and $\pi:A\rightarrow A/I$ is the qutiont map.$\hfill\square$

\section{\normalsize\bf Cohomological properties of projective tensor products for Banach algebras}

For Banach algebras $A_1$ and $A_2$, let $A_1\hat{\otimes}A_2$  be the projective tensor product.

\begin{theorem}\label{ala}  Let $A_1$ and $A_2$ be  Banach algebras. Then $A_1$ and $A_2$ are 0-point amenable if and only if  $A_1\hat{\otimes}A_2$ is 0-point amenable.
\end{theorem}
{\it Proof.}  Assume first that $A_1$ and $A_2$ are unital. Let $d: A_1\hat{\otimes}A_2\rightarrow\Bbb{C}$ be a continuous point derivation at $\varphi\in\Delta_0(A_1\hat{\otimes}A_2)$ of $A_1\hat{\otimes}A_2$. Define $d_1: A_1\rightarrow{\Bbb C}$ by $$d_1(a)=d(a\otimes 1).$$ Then $d_1$ is a continuous point derivation at $\varphi_1\in\Delta_0(A_1)$ of $A_1$, where $\varphi_1$ is defined by $\varphi_1(a_1)=\varphi(a_1\otimes 1)$ for all $a_1\in A_1$. Similarly, $d_2$ defined by $d_2(a_2)=d(1\otimes a_2)$ is a continuous point derivation  at $\varphi_2\in\Delta_0(A_2)$ of $A_2$, where $\varphi_2(a_2)=\varphi(1\otimes a_2)$ for all $a_2\in A_2$. If $A_1$ and $A_2$ are 0-point amenable, then $d_1$ and $d_2$ are zero. Hence for every $a_1\in A_1$ and $a_2\in A_2$ we have
\begin{eqnarray*}
d(a_1\otimes a_2)&=&d((a_1\otimes 1)(1\otimes a_2))\\
&=& d(a_1\otimes 1)\varphi(1\otimes a_2)+\varphi(a_1\otimes 1)d(1\otimes a_2)\\
&=&d_1(a_1)\varphi_2(a_2)+\varphi_1(a_1)d_2(a_2)\\
&=&0.
\end{eqnarray*}
Hence $A_1\hat{\otimes}A_2$ is 0-point amenable.

Now, let $A_1$ and $A_2$ be any Banach algebras. Assume that  $A_1$ and $A_2$ are 0-point amenable. It follows from Theorem \ref{sandali124} that $A_1^\sharp$ and $A_2^\sharp$ are 0-point amenable and so  $A_1^\sharp\hat{\otimes}A_2^\sharp$ is  0-point amenable. But,
\begin{eqnarray}\label{tem}
A_1^\sharp\hat{\otimes}A_2^\sharp&=&(A_1\oplus{\Bbb C})\hat{\otimes}(A_2\oplus{\Bbb C})\nonumber\\
&=&(A_1\hat{\otimes}A_2)\oplus(A_1\hat{\otimes}{\Bbb C})\\
&\oplus&({\Bbb C}\hat{\otimes}A_2)\oplus({\Bbb C}\hat{\otimes}{\Bbb C}).\nonumber
\end{eqnarray}
This together with Theorem \ref{sandali124}, $A_1\hat{\otimes}A_2$ is 0-point amenable.

To prove the converse, let $d_1: A_1\rightarrow \Bbb{C}$ be a continuous point derivation at $\varphi_1\in\Delta_0(A_1)$. Define $d: A_1\hat{\otimes}A_2\rightarrow\Bbb{C}$ by $d(a_1\otimes a_2)=d_1(a_1)\varphi_2(a_2)$, where $\varphi_2\in\Delta_0(A_2)$. Then $d$ is a countinuous point derivation at $\varphi\in\Delta_0(A_1\hat{\otimes}A_2)$ of $A_1\hat{\otimes}A_2$, where
$\varphi(a_1\otimes a_2)=\varphi_1(a_1)\varphi_2(a_2)$ for all $a_1\in A_1$ and $a_2\in A_2$. Since $A_1\hat{\otimes}A_2$ is 0-point amenable, $d=0$ and hence $d_1=0$. Thus $A_1$ is 0-point amenable. Similarly, $A_2$ is 0-point amenable. $\hfill\square$\\

As an immediate consequence of Theorem 3.3 in \cite{mr4} and Theorem \ref{ala} we have the following result.

\begin{corollary}\label{sal8} Let $A_i$  be a Banach algebra with  $\Delta(A_i)$ is a non-empty set, for $i=1, 2$. Then the following assertions are equivalent.

\emph{(a)} $A_1\hat{\otimes}A_2$ is cyclically weakly amenable.

\emph{(b)}  $A_1$ and $A_2$ are cyclically weakly amenable.

\emph{(c)} $A_1^\sharp$ and $A_2^\sharp$ are cyclically weakly  amenable.

\emph{(d)} $A_1^\sharp\hat{\otimes}A_2^\sharp$ is cyclically weakly amenable.
\end{corollary}

Groenbaek \cite{gro} proved some hereditary properties of weak amenability for commutative Banach algebras. He also showed that if $A_1$ and $A_2$ are commutative weakly amenable Banach algebra, then $A_1\hat{\otimes} A_2$ is weakly amenable. The converse of this result is proved by  Yazdanpanah \cite{y}. He also found some result concerning weak amenability of tensor product of noncommutative Banach algebras.

\begin{theorem}\label{end} Let $A_i$ be a unital Banach algebra with $\Delta(A_i)\neq\emptyset$, for $i=1, 2$. Then the following statements hold.

\emph{(i)} If $A_1\hat{\otimes}A_2$ is cyclically amenable, then $A_1$ and $A_2$ are cyclically amenable.

\emph{(ii)} If $A_1\hat{\otimes}A_2$ is weakly amenable, then $A_1$ and $A_2$ are weakly amenable.
\end{theorem}
{\it Proof.} Choose $\varphi_2\in\Delta(A_2)$ and define the function $\Lambda:A_1\hat{\otimes}A_2\rightarrow A_1$ by $$\Lambda(a_1\otimes a_2)=\varphi_2(a_2)a_1.$$ Then $\Lambda$ is an epimorphism. If we define the function $\Gamma: A_1\rightarrow A_1\hat{\otimes}A_2$ by $$\Gamma(a_1)=a_1\otimes 1,$$ then $\Lambda\circ\Gamma=id_{A_1}$, the identity map on $A_1$. Thus $\Lambda$ is a retraction.  Now, apply Theorem \ref{react} to conclude that cyclic amenability of  $A_1\hat{\otimes}A_2$ forces  $A_1$ to be cyclically amenable. Similarly, $A_2$ is weakly amenable. So (i) holds. The statements (i) together with Corollary \ref{sal8} proves (ii).$\hfill\square$

\section{\normalsize\bf Type of weak amenability of the second dual
Banach algebras}

For a Banach algebra $A$, Gourdeau \cite{g} proved that if $A^{**}$ is amenable, then so dose $A$; see \cite{glw} for another proof of it. Ghahramani, Loy and Willis \cite{glw} investigated this result for weak amenability instead of amenability. They showed that if $A$ is a left ideal in $A^{**}$, then Gourdea's theorem remains valid for weak amenability instead of amenability. They claimed that if every derivation from $A$ into $A^*$ is weakly compact, then weak amenability of $A^{**}$ implies weak amenability of $A$. Several authors doubted the correctness of it and proved that under certain circumstances the claim is true \cite{ap, aef, drv, ef, gl}. For example, Dales, Rodriguez-Palacios and Velasco \cite{drv} gave conditions under which the second transpose of a continuous derivation $D: A\rightarrow A^*$ is also a derivation. Then, they proved the claime under the additional assumption that $A$ is Arens regular. Finally, Ghahramani, Loy and Willis \cite{glw1} proved it for any Banach algebra; see also \cite{bv}. In the following, we prove this result for cyclic and point  amenability.

\begin{proposition} Let $A$ be a Banach algebra and every continuous derivation from a Banach algebra $A$ into $A^*$ is weakly compact. Then the following statements hold.

\emph{(i)} A bounded linear operator $D: A\rightarrow A^*$ is a cyclic derivation if and only if $D^{**}$ is a cyclic derivation.

\emph{(ii)} If $A^{**}$ is  cyclically weakly amenable, then $A$ is  cyclically weakly amenable.

\emph{(iii)} Assume that $A$ is a unital, commutative Banach algebra. If $A^{**}$ is cyclically amenable, then $A$ is cyclically amenable.
\end{proposition}
{\it Proof.} Let $D: A\rightarrow A^*$ be a continuous derivation. By assumption, $D$ is weakly compact and so $D^{**}$ is a continuous derivation. If $\Phi, \Psi\in A^{**}$, then there exist nets $(a_\alpha)$ and $(b_\beta)$ in $A$ such that $a_\alpha\rightarrow\Phi$ and $b_\beta\rightarrow \Psi$ in the weak$^*$ topology of $A^{**}$. Hence
\begin{eqnarray}\label{bbb}
\lim_\alpha\lim_\beta\langle D(a_\alpha), b_\beta\rangle=\langle D^{**}(\Phi), \Psi\rangle.
\end{eqnarray}
Assume that $D$ is cyclic. Then (\ref{bbb}) implies that $D^{**}$ is cyclic. The converse follows from the fact that $A$ is a subspace of $A^{**}$. So (i) holds. The statement (ii) follows at once from (i). Finally, let $D: A\rightarrow A^*$ be a continuous cyclic derivation. Then for every $a\in A$, we have
$$
\langle D^*(1_A), a\rangle=\langle D(a), 1_A\rangle=0
$$
and so for every $\Phi\in A^{**}$, we have
$$
\langle D^{**}(\Phi), 1_A\rangle=\langle\Phi, D^*(1_A)\rangle=0.
$$
Hence $D^{**}$ is cyclic. By assumption $D^{**}=0$. Therefore, $D=0$.
$\hfill\square$

\begin{proposition}\label{pp} Let $A$ be a Banach algebra. Then the following statements hold.

\emph{(i)} A bounded linear  functional $d\in A^*$  is a continuous point derivation of $A$ at $\varphi\in\Delta(A)$ if and only if  $d^{**}$  is a continuous point derivation of $A^{**}$ at $\varphi^{**}\in\Delta(A^{**})$

\emph{(ii)} If $A^{**}$ is point amenable, then $A$ is point amenable.

\emph{(ii)} If $A^{**}$ is 0-point amenable, then $A$ is 0-point amenable.
\end{proposition}
{\it Proof.} This is routine.$\hfill\square$

\begin{corollary} Let the second dual of Banach algebra $A$ be weakly amenable. Then the following statements hold.

\emph{(i)} $A$ is 0-point amenable.

\emph{(ii)} If $\Delta(A)\neq\emptyset$, then $A$ is cyclically weakly amenable.

\emph{(iii)} If $A$ is cyclic amenable and $\Delta(A)\neq\emptyset$, then $A$ is weakly amenable.

\emph{(iv)} If $A$ is unital and commutative, then $A$ is weakly amenable.

\emph{(v)} If $A$ is semisimple, then $A$ is weakly amenable.
\end{corollary}
{\it Proof.} Let $A^{**}$ be weakly amenable. Then $A^{**}$ is 0-point amenable; see Theorem 4.1 in \cite{mr4}. In view of Proposition \ref{pp}, $A$ is 0-point amenable. So (i) holds. For (ii) and (iii), we only recall from \cite{mr4} that if $\Delta(A)$ is non-empty, then cyclically weak amenability and point amenability are equivalent. Also, $A$ is weakly amenable if and only if $A$ is cyclically amenable and cyclically weakly amenable. Finally, the statements (iv) and (v) follow from the fact that  if $A$ is either a unital commutative Banach algebra or semisimple, then weak amenability and point amenability coincide \cite{mr4}.$\hfill\square$\\

Let $A$ be a Banach algebra with the multiplication ``$\cdot$''. Then $A$ with the norm $\|.\|_A$ and the multiplication ``$\bullet$'' defined by
$$
a\bullet b=b\cdot a\quad\quad(a, b\in A)
$$
is also a Banach algebra. This Banach algebra is called the \emph{opposite algebra} \emph{of} $A$ and is denoted by $A^{op}$.

\begin{theorem}\label{op} Let $A$ be a Banach algebra. Then $A$ is weakly (respectively, point, cyclically) amenable if and only $A^{op}$ is weakly (respectively, point, cyclically) amenable.
\end{theorem}
{\it Proof.} Let $D: A^{op}\rightarrow (A^{op})^*$ be a continuous derivation. We define the bounded linear operator $\tilde{D}: A\rightarrow A^*$ by $$
\langle \tilde{D}(a), b\rangle=\langle D(a), b\rangle\quad\quad(a,b\in A).
$$
Then for every $a, b, c\in A$ we have
\begin{eqnarray*}
\langle \tilde{D} (a\cdot b), c\rangle&=&\langle D(b\bullet a), c\rangle\\
&=&\langle b\bullet D(a), c\rangle+\langle D(b)\bullet a, c\rangle\\
&=&\langle D( a), c\bullet b\rangle+ \langle D(b), a\bullet c\rangle\\
&=&\langle D(a), b\cdot c\rangle+\langle D(b), c\cdot a\rangle\\
&=&\langle \tilde{D} (a), b\cdot c\rangle+\langle \tilde{D} (b), c\cdot a\rangle\\
&=&\langle \tilde{D} (a)\cdot b, c\rangle+\langle a\cdot \tilde{D} (b), c\rangle.
\end{eqnarray*}
Hence $\tilde{D}$ is a continuous derivation. If $A$ is weakly amenable, then
$D=\hbox{ad}_F$ for some $F\in A^*$. Thus
\begin{eqnarray*}
\langle D(a), b\rangle&=&
F(a\bullet b)-F(b\bullet a)\\
&=&F(b\cdot a)-F(a\cdot b)=\langle\hbox{ad}_{-F}(a), b\rangle.
\end{eqnarray*}
It follows that $D$ is inner. Therefore, $A^{op}$ is weakly amenable. The other case is proved similarly. For the converse, we only note that $(A^{op})^{op}=A$.$\hfill\square$

\vspace{2mm}

 {\footnotesize
\noindent {\bf Mohammad Javad Mehdipour}\\
Department of Mathematics,\\ Shiraz University of Technology,\\
Shiraz
71555-313, Iran\\ e-mail: mehdipour@.ac.ir\\
{\bf Ali Rejali}\\
Department of Pure Mathematics,\\ Faculty of Mathematics and Statistics,\\ University of Isfahan,\\
Isfahan
81746-73441, Iran\\ e-mail: rejali@sci.ui.ac.ir\\

\footnotesize

\begin{thebibliography}{99}

\bibitem{ap} M. J. Aleandro and C. C. Pena, On dual-valued operators on Banach algebras, New York J. Math., 18 (2012) 657--665.

\bibitem{aef} M. Amini, M. Essmaili and M. Filali, The second transpose of a derivation and weak amenability of the second dual Banach algebras, New York J. Math., 22 (2016) 265--275

\bibitem{bcd} W. Bade, P. C. Curtis and H. G. Dales, Amenability and weak amenability for Beurling and Lipschitz algebras, Proc. London Math. Soc., (3) 55 (1987) 359--377.

\bibitem{bv} S. Barootkoob and H. R. Ebrahimi Vishki, Lifting derivations and n-weak amenability of
the second dual of a Banach algebra, Bull. Aust. Math. Soc., 83 (2011) 122–129.

\bibitem{cl} P. C. Curtis and R. J. Loy, The structure of amenable Banach algebras, J. London Math. Soc., 40 (1) (1989) 89--104.

\bibitem{d} H. G. Dales, Banach algebras and Automatic Continuity, Clarendon Press, Oxford, 2000.

\bibitem{drv} H. G. Dales, A. Rodriguez-Palacios and M. V. Velasco, The second transpose of a derivation, J. London Math. Soc., 64 (3) (2001) 707--721.

\bibitem{ef} M. Eshaghi Gordji and M. Filali, Weak amenability of the second dual of a Banach algebra, Studia Math., 182 (3) (2007) 205--213.

\bibitem{f} B. Forrest, Weak amenability and the second dual of the Fourier algebra, Proc. Amer. Math. Soc., 125 (8) (1997) 2373--2378.

\bibitem{gl} F. Ghahramani and J. Laali, Amenability and topological centres of the second duals of Banach algebras, Bull. Austral. Math. Soc., 65 (2) (2002) 191--197.

\bibitem{glw} F. Ghahramani, R. J. Loy and G. A. Willis, Amenability and weak amenability of second conjugate Banach algebras, Proc. Amer. Math. Soc., 124 (5)  (1996) 1489--1497.

\bibitem{glw1} F. Ghahramani, R. J. Loy and G. A. Willis, Addendum to "Amenability and weak amenability of second conjugate Banach algebras'', Proc. Amer. Math. Soc., 148 (10) (2020) 4573--4575.

\bibitem{g} F. Gourdeau, Amenability and the second dual of a Banach algebra, Studia Math., 125 (1) (1997) 75--81.

\bibitem{gro} N. Gronback, A charactrization of weakly amenable Banach algebras, Studia. Math., 94 (1989) 149--162.

\bibitem{gro1} N. Gronback, Amenability and weak amenability of tensor algebras and algebras of nuclear operators, J. Austral. Math. Soc., (Series A) 51 (1991) 483--488.

\bibitem{gro2} N. Gronback, Weak and cyclic amenability for non-commutative Banach algebras, Proc.
Edinburgh Math. Soc.,  35 (1992) 315--328.

\bibitem{j1} B. E. Johnson, Derivations from $L^1(G)$ into $L^1(G)$ and $L^\infty(G)$, Harmonic analysis (Luxembourg, 1987), 191--198, Lecture Notes in Math., 1359, Springer, Berlin, 1988.

\bibitem{ll} A. T. Lau and R. J.  Loy,  Weak amenability of Banach algebras on locally compact groups, J. Funct. Anal., 145 (1) (1997) 175--204.

\bibitem{mr4} M. J. Mehdipour and A. Rejali, Different types of weak amenability for Banach algebras,  arXiv:2209.13580.

\bibitem{m} T. O. Mewomo,  Various notions of amenability in Banach algebras, Expo. Math., 29 (3) (2011) 283--299.

\bibitem{pat} A. L. T. Paterson, Virtual diagonals and $n-$amenability for Banach algebras, Pacific
J. Math., 175 (1996) 161–-185.

\bibitem{mon} M. Sanjani Monfared, On certain products of Banach algebras with
applications to harmonic analysis, Studia  Math., 178 (3) (2007) 277--294.

\bibitem{y} T. Yazdanpanah, Weak amenability of tensor product of Banach algebras, Proc. Rom. Acad. Ser. A Math. Phys. Tech. Sci. Inf. Sci., 13 (4) (2012) 310--313.
\end{thebibliography}
\end{document}